\documentclass{ieeeconf}
\IEEEoverridecommandlockouts
\usepackage{verbatim}
\usepackage{epsfig}
\usepackage{amsmath}
\usepackage{amsthm}
\usepackage{amssymb}
\usepackage{psfrag}
\usepackage[T1]{fontenc}
\usepackage[ruled]{algorithm2e}

\usepackage{url}
\usepackage{dsfont}
\usepackage{mathtools}
\DeclarePairedDelimiter{\ceil}{\lceil}{\rceil}

\usepackage[usenames,dvipsnames]{pstricks}
\usepackage{pst-grad} 
\usepackage{pst-plot} 
\usepackage[ruled]{algorithm2e}
\usepackage{etoolbox}

\usepackage{graphicx}

\usepackage[font=footnotesize,skip=0pt]{caption}

\makeatletter
\patchcmd{\@makecaption}
  {\scshape}
  {}
  {}
  {}
\makeatother

\usepackage[noadjust]{cite}

\newtheorem{theorem}{Theorem}
\newtheorem{lemma}{Lemma}

\newtheorem{definition}{Definition}
\newtheorem{assmp}{Assumption}

\newtheorem{proposition}{Proposition}

\newcommand{\mfi}{M_{\mu f_i}}
\newcommand{\mgi}{M_{\mu g_i}}
\newcommand{\mFi}{F_{i,\mu}}

\newcommand{\mF}{\frac{1}{n} \sum_{i=1}^n F_{i,\mu}}
\newcommand{\y}{\widehat{y}}

\newcommand{\gmF}{\frac{1}{n} \sum_{i=1}^n \nabla F_{i,\mu}}

\newcommand{\dc}{\texttt{DDC-Consensus}}
\newcommand{\dci}{\texttt{DDC-Consensus-Inexact}}
\newcommand{\dcm}{\texttt{DDC-Mixing}}

\newcommand{\V}{\mathcal{V}}

\DeclareMathOperator*{\argmin}{arg\ minimize}
\DeclareMathOperator*{\minimize}{minimize}

\DeclareMathOperator*{\dom}{dom}

\IEEEpeerreviewmaketitle

\begin{document}

\title{\textbf{Distributed Difference of Convex Optimization}

}
\author{Vivek Khatana, Murti V. Salapaka
\thanks{ Vivek Khatana (Email: \{khata010\}@umn.edu) and Murti V. Salapaka (Email: \{murtis\}@umn.edu) are with the Department of Electrical and Computer Engineering, University of Minnesota, Minneapolis, USA.}
}

\maketitle

\begin{abstract}
In this article, we focus on solving a class of distributed optimization problems involving $n$ agents with the local objective function at every agent $i$ given by the difference of two convex functions $f_i$ and $g_i$ (difference-of-convex (DC) form), where $f_i$ and $g_i$ are potentially nonsmooth. The agents communicate via a directed graph containing $n$ nodes. We create smooth approximations of the functions $f_i$ and $g_i$ and develop a distributed algorithm utilizing the gradients of the smooth surrogates and a finite-time approximate consensus protocol. We term this algorithm as $\dc$. The developed $\dc$ algorithm allows for non-symmetric directed graph topologies and can be synthesized distributively. We establish that the $\dc$ algorithm converges to a stationary point of the nonconvex distributed optimization problem. The performance of the $\dc$ algorithm is evaluated via a simulation study to solve a nonconvex DC-regularized distributed least squares problem. The numerical results corroborate the efficacy of the proposed algorithm.\\[1ex]
\textit{keywords}: Distributed optimization, nonconvex optimization, difference-of-convex (DC) functions, DC programming, distributed gradient descent, directed graphs.
\end{abstract}

\section{Introduction}
Due to the increase in size and complexity of modern systems, solving problems involving many agents via distributed methods is highly desirable. An effective way towards a distributed solution is to cast the problems in the framework of distributed optimization \cite{ddistadmm_cdc, gradcons_acc}. In this article, we focus on the following distributed optimization problem,
\begin{align}\label{eq:introprob}
    \textstyle \minimize \limits_{x \in \mathbb{R}^p} \ F(x) = \frac{1}{n}\sum_{i=1}^{n} ( f_i(x) - g_i(x) ),
\end{align}
where  $x \in \mathbb{R}^p$ is a common decision, and $f_i: \mathbb{R}^p \to \mathbb{R} \cup \{-\infty, \infty\}$ and $g_i: \mathbb{R}^p \to \mathbb{R}$ are private objective functions of agent $i$. The agents are connected through a directed graph, $\mathcal{G}(\mathcal{V},\mathcal{E})$,  where $\mathcal{V}$ and $\mathcal{E}$ are the set of vertices and edges respectively. The terms in the aggregate objective function in~\eqref{eq:introprob} are of a difference-of-convex (DC) form. Due to the richness of the set of DC functions (see the properties mentioned in \cite{hartman1959functions}), DC functions can be used to model most of all practical non-convex
optimization problems. Many applications in statistical learning and estimation \cite{nouiehed2019pervasiveness}, power systems \cite{merkli2017fast}, computational biology \cite{an2003solving}, signal restoration \cite{an2001dc}, network optimization \cite{hoai2002dc}, combinatorial optimization \cite{le2001continuous} can be posed as DC programs.

The initial study of distributed optimization can be traced back to the seminal works \cite{tsitsiklis1984problems, bertsekas1989parallel}. Since then the special case with $f_i$ being convex and $g_i = 0$ has received significant development (see \cite{pu2020push, khatana_admm_tac, gradcons_tnse} and references therein). As a relaxation of the convex objective functions most works in the literature are restricted to Lipschitz differentiable not necessarily convex functions $f_i$ with $g_i = 0$ (some examples include \cite{non_cvx_wotao_yin, Xin2020AFR, Yan2020DistributedNO, scutari2019distributed, gorbunov2021marina, Jie_prox, Lorenzo2016NEXTIN}) and do not address the general problem~\eqref{eq:introprob}. Article \cite{Chen2020OnDN} focuses on a special case with weakly convex objective functions $f_i$ and $g_i = 0$. Article \cite{alvarado2014new} considers a DC form however the functions $f_i$ and $g_i$ are assumed to have Lipschitz differentiable gradients, and does not ensure globally agreed decisions among all the agents. 
To the best of the authors' knowledge, the general
\textit{distributed} optimization problem~(\ref{eq:introprob}) with DC functions is not addressed in prior works.

However, DC programming is studied in the centralized optimization literature \cite{dchomepage, le2015dc, tao1997convex, an2005dc, ahn2017difference, nhat2018accelerated, sun2018alternating}. Articles \cite{dchomepage, le2015dc, tao1997convex, an2005dc, ahn2017difference} utilize an affine minorization $g_i^m$ of the functions $g_i, i = 1$ and minimize the resulting convex function $f_i - g_i^m, i = 1$. The work in \cite{nhat2018accelerated} provides an accelerated version of the algorithms proposed in \cite{dchomepage, le2015dc, tao1997convex, an2005dc, ahn2017difference} using adaptive step-sizes and extrapolation of the algorithm iterates. Article \cite{sun2018alternating} proposed an alternating direction of multipliers method based algorithm for solving centralized DC optimization problems with Lipschitz differentiable functions and a linear equality constraint. The current article is closest to the centralized algorithms \cite{sun2023algorithms, wen2018proximal}, where a majorization of the objective functions via the proximal mapping is solved.

In this article, we propose an iterative algorithm termed \underline{\textbf{D}}istributed \underline{\textbf{DC}}-Consensus $\dc$. The $\dc$ algorithm proceeds in two steps: first, all the agents in parallel solve a local optimization problem involving $f_i - g_i$, then each agent $i \in \V$ shares the obtained solution with its neighbors where the estimates of the global solution are updated by utilizing an average consensus protocol (the algorithm is presented in detail in Section~\ref{sec:proposed_alg}). The main contributions of the current article are as follows: \\[0.1ex]
\textbf{i)} We develop a distributed algorithm, $\dc$, based on distributed gradient descent to the solve~\eqref{eq:introprob} for non-differentiable weakly convex $f_i$ and non-differentiable convex function $g_i$. These assumptions are weaker than the existing literature on distributed optimization.\\[0.1ex]
\textbf{ii)} We establish that the $\dc$ converges to a stationary solution of the nonconvex problem~(\ref{eq:introprob}).\\[0.1ex]
\textbf{iii)} The proposed $\dc$ algorithm has desirable properties as it (a) is suitable for directed graphs and non-symmetric communication topologies, and (b) can be implemented as well as synthesized distributively.\\[0.1ex]
We demonstrate the performance of the $\dc$ algorithm by solving a 
DC-regularized least squares problem. Empirical results from the numerical simulation study corroborate that the proposed $\dc$ algorithm performs well in practice.

The rest of the article is organized as follows: Some essential definitions and notations are presented in Section~\ref{sec:defn}. Section~\ref{sec:prob_form} provides the details of the proposed $\dc$ algorithm. Convergence analysis of $\dc$ is presented in Section~\ref{sec:convergence}. Section~\ref{sec:simulation} presents the numerical simulation study followed by concluding remarks in Section~\ref{sec:conclusion}.

\vspace{-0.1in}
\subsection{Definitions and Notations}\label{sec:defn}
\begin{definition}(Directed Graph)
A directed graph $\mathcal{G}$ is a pair $(\mathcal{V},\mathcal{E})$ where $\mathcal{V}$ is a set of vertices and $\mathcal{E}$ is a set of edges, which are ordered subsets of two distinct elements of $\mathcal{V}$. If an edge from $j \in \mathcal{V}$ to $i \in \mathcal{V}$ exists it is denoted as $(i,j)\in \mathcal{E}$. 
\end{definition}

\begin{definition}(Strongly Connected Graph) A directed graph is strongly connected if for any pair $(i,j),\ i\not =j$, there is a directed path from node $i$ to node $j$.
\end{definition}

\begin{definition}(Diameter of a Graph)
The longest shortest directed path between any two nodes in the graph. 
\end{definition}

\begin{definition}(In-Neighborhood) The set of in-neighbors of node $i \in \mathcal{V}$ is called the in-neighborhood of node $i$ and is denoted by $\mathcal{N}^{-}_i = \{j \ | \ (i,j)\in \mathcal{E}\}$ not including the node $i$.
\end{definition}

\begin{definition}(Column-Stochastic Matrix) A matrix $M=[m_{ij}] \in \mathbb{R}^{n\times n}$ 
is called a column-stochastic matrix if $0 \leq m_{ij}\leq 1$ and $\sum_{i=1}^{n}m_{ij}=1$ for all $ 1 \leq i,j \leq n$. 
\end{definition}


\begin{definition}(Lipschitz Continuity and Differentiability)
A function $f: \mathbb{R}^p \to \mathbb{R}$ is called Lipschitz continuous with constant $L > 0$ and Lipschitz differentiable with constant $L_f > 0$ respectively, if the following inequalities hold: 
\begin{align*}
    | f(x) - f(y) | &\leq L \|x - y\|, \ \forall \ x,y \in \mathbb{R}^p \\
    \| \nabla f(x) - \nabla f(y)\| &\leq L_f \|x-y\|, \ \forall \ x,y \in \mathbb{R}^p.
\end{align*}
\end{definition}

\begin{definition}(Strongly Convex Function) A differentiable function $f$ is called strongly convex with parameter $\sigma$, if there exists $\sigma > 0$ such that for all $x, y$ in the domain of $f$:
\begin{align*}
   \langle \nabla f(x)-\nabla f(y), x-y \rangle \geq \sigma \|x-y\|^2, \ \forall \ x,y \in \dom f.
\end{align*}
\end{definition}

\begin{definition}(Weakly Convex Function) For some $m_f \geq 0$, we say that function $f: \mathbb{R}^p \to \mathbb{R}$ is $m_f$-weakly convex if $ f + \frac{m_f}{2} \|x\|^2$ is convex. 
\end{definition}

\begin{definition}(Level-bounded Function) Function $f: \mathbb{R}^p \to \mathbb{R}$ is called level-bounded if the $\beta$-level-set $\{ x : f(x) \leq \beta \}$ is bounded (possibly empty) for all $\beta \in \mathbb{R}$.
\end{definition}

\begin{definition}(Lower Semicontinuity) A function $f: \mathbb{R}^p \to \mathbb{R} \cup \{-\infty, \infty\}$  is called lower semi-continuous at a point $x_0$ if for every real $y < f(x_0)$ there exists a neighborhood $U$ of $x_0$ such that $y < f(x)$ for all $x \in U$. The function $f$ is lower semicontinuous at every point of its domain, $\dom(f)$.
\end{definition}

\begin{definition}(Proper Function) A function $f: \mathbb{R}^p \to \mathbb{R} \cup \{-\infty, \infty\}$  is called proper if $f(x) > -\infty$ for every $x$ and and if there also exists \textit{some} point $x_0$ such that $f(x_0) < \infty$.
\end{definition}

\noindent $\|x\|_1$ and $\|x\|$ denote the 1-norm and 2-norm of the vector $x \in \mathbb{R}^p$ respectively. The notation $\ceil*.$ denotes the least integer function or the ceiling function, defined as:
given $x \in \mathbb{R}, \ceil*x = \min \{ m \in \mathbb{Z} | m \geq x\},$ where $\mathbb{Z}$ is the set of integers. $\dom (f) = \{ x: f(x) < \infty \}$.

\section{Proposed Methodology}\label{sec:prob_form}
\noindent The following assumption is satisfied throughout the article,
\begin{assmp}\label{assm:weakly_convex}
    1. Functions $f_i : \mathbb{R}^p \to \mathbb{R} \cup \{-\infty, \infty \}$ are proper, lower semi-continuous, and $m_{f_i}$-weakly convex. \\
    2. Functions $g_i : \mathbb{R}^n \to \mathbb{R}$ are convex and  finite everywhere.\\
    3. The set of global minimizers of~\eqref{eq:introprob}, $\argmin_x F(x)$, is nonempty, and the global minimum value $F^*$ is finite.
\end{assmp}
\noindent A vector $\tilde{y} \in \mathbb{R}^p$ is called a stationary point of $F$ if:
\begin{align}
   \textstyle 0 \in \partial \left( \frac{1}{n} \sum_{i=1}^n f_i(\tilde{y}) \right) - \partial \left( \frac{1}{n} \sum_{i=1}^n g_i(\tilde{y}) \right), \label{eq:stationary_condition}
\end{align}
or equivalently, $\partial \frac{1}{n} \sum_{i=1}^n f_i(\tilde{y}) \cap \partial \frac{1}{n} \sum_{i=1}^n g_i(\tilde{y}) \neq 0$. Furthermore, given $ \varepsilon > 0$, we say $\tilde{y} \in \mathbb{R}^p$ is an
$ \varepsilon$-stationary point of $F$ if there exists $(\xi; y) \in \mathbb{R}^p \times \mathbb{R}^p$ such that
\begin{align}
    \textstyle \xi \in \partial \left( \frac{1}{n} \sum_{i=1}^n f_i(\tilde{y}) \right) - \ & \textstyle \partial \left( \frac{1}{n} \sum_{i=1}^n g_i(y) \right),  \ \mbox{and} \nonumber \\
    \max \big\{ \|\xi\|, & \| \tilde{y} - y\| \big \} \leq \varepsilon. \label{eq:epsilon_stationary_condition}
\end{align}
We use $\partial f(x)$ to denote the general subdifferential of the function $f$ at $x$ (\cite{rockafellar2009variational}, Definition $8.3$). 

\subsection{Approach Roadmap and Supporting Results}
Note that the objective function has a DC structure. Under Assumption~\ref{assm:weakly_convex} we do not impose differentiability of the functions $f_i$ and $g_i$. Thus, to develop an algorithm with desirable convergence properties, we first utilize a smooth function mapping to obtain a differentiable DC approximation of the objective function $f_i - g_i$. Then we utilize the differentiable DC approximation to develop first-order algorithms to solve problem~\eqref{eq:introprob}. We take the approach in \cite{sun2023algorithms} and utilize separate Moreau envelopes of $f_i$ and $g_i$ to obtain a smooth approximation of the function $f_i - g_i$. To this end, we define the Moreau envelops: Given $ 0 < \mu_i < 1/m_{f_i}$, the Moreau envelopes, $\mfi: \mathbb{R}^p \to \mathbb{R}$, and $\mgi: \mathbb{R}^p \to \mathbb{R}$ of functions $f_i$ and $g_i$ respectively are given by    
\begin{align}
   M_{\mu_i f_i}(y) &:= \minimize_{x \in \mathbb{R}^p} \textstyle \left \{ f_i(x) + \frac{1}{2 \mu_i} \| x - y \|^2 \right \}, \label{eq:individual_moreau_f} \\
   M_{\mu_i g_i}(y) &:=  \minimize_{x \in \mathbb{R}^p} \textstyle \left \{ g_i(x) + \frac{1}{2 \mu_i} \| x - y \|^2 \right \},\label{eq:individual_moreau_g}
\end{align}
and form the smooth function $F_{i, \mu_i} := M_{\mu_i f_i} - M_{\mu_i g_i}$, for all $i \in \{ 1,\dots,n \}$. The corresponding proximal mappings 
$x_{\mu_i f_i}: \mathbb{R}^p \to \mathbb{R}^p$, and $x_{\mu_i g_i}: \mathbb{R}^p \to \mathbb{R}^p$ are defined as,
\begin{align}
   x_{\mu_i f_i}(y) &:= \argmin_{x \in \mathbb{R}^p} \textstyle  \left \{ f_i(x) + \frac{1}{2 \mu_i} \| x - y \|^2 \right \}, \label{eq:individual_prox_maps_f} \\
   x_{\mu_i g_i}(y) &:= \argmin_{x \in \mathbb{R}^p} \textstyle  \left \{ g_i(x) + \frac{1}{2 \mu_i} \| x - y \|^2 \right \}. \label{eq:individual_prox_maps_g}
\end{align}

\noindent Define, $m_f := \max_{1 \leq i \leq n} m_{f_i}$. Given, $ 0 < \mu < 1/m_f$, let $F_{i,\mu}(x) := \mfi(x) - \mgi(x)$, and 
\begin{align}\label{eq:f_mu}
    \hspace{-0.1in} F_{\mu}(x) := \mF(x) = \frac{1}{n} \sum_{i=1}^n (\mfi(x) - \mgi(x)).
\end{align}
Next, we summarize the properties of $F_{i,\mu}$, the Moreau envelops $\mfi, \mgi$ the mappings $x_{\mu f_i}, x_{\mu g_i}$. 
 
\begin{lemma} \label{lem:moreau_properties} (\cite{sun2023algorithms}, Properties of Moreau envelops and proximal mappings)  
Let Assumption~\ref{assm:weakly_convex} holds. Let $0 < \mu < 1/m_f$ and $\mfi, \mgi, x_{\mu f_i}, x_{\mu g_i}$ be given by definitions~\eqref{eq:individual_moreau_f}-\eqref{eq:individual_prox_maps_g} with parameter $\mu$. Then, the following claims hold:\\
1. $x_{\mu f_i}, x_{\mu g_i}$ are Lipschitz continuous with modulus, $\frac{1}{1 - \mu m_f}$ and $1$ respectively. \\
2. $\mfi, \mgi$ are differentiable with gradient $\nabla \mfi(y) = \mu^{-1} ( y - x_{\mu f_i}(y) ), \nabla \mgi(y) = \mu^{-1} ( y - x_{\mu g_i}(y) )$.\\
3. $\nabla \mfi, \nabla \mgi$ are Lipschitz continuous with modulus $L_{\mu_{M_f}} = \frac{2 - \mu m_f}{\mu - \mu^2 m_f}$ and $ L_{\mu_{M_g}} = \frac{2}{\mu}$ respectively. \\
4. $F_{i,\mu}$ is differentiable, and $\nabla F_{i,\mu}(y) = \nabla \mfi(y) - \nabla \mgi(y) = \mu^{-1} (x_{\mu g_i}(y) - x_{\mu f_i}(y) )$ is Lipschitz continuous with the modulus $L_{\mu_F} := \frac{2 - \mu m_f}{\mu - \mu^2 m_f}$.
\end{lemma}

\noindent Next, we introduce the following optimization problem,
\begin{align}\label{eq:intro_mu_prob}
    \textstyle \minimize \limits_{x \in \mathbb{R}^p} \ F_{\mu}(x) = \frac{1}{n} \sum_{i=1}^n (\mfi(x) - \mgi(x)). 
\end{align}

\begin{lemma}\label{lem:problem_correspondance} (\cite{sun2023algorithms}, Relation between~\eqref{eq:introprob} and~\eqref{eq:intro_mu_prob}) Let Assumption~\ref{assm:weakly_convex} holds and $0 < \mu < 1/m_f$. Then,\\
1. The set of global minimizers of~\eqref{eq:intro_mu_prob}, $\argmin_{x} F_\mu(x)$, is non-empty, and $F^* = \minimize_{x \in \mathbb{R}^p} F_\mu(x)$. \\
2. $\tilde{y}$ is a stationary point of $F_\mu$, i.e. $\nabla F_\mu(\tilde{y}) = 0$, if and only if $ \tilde{x} = \frac{1}{n} \sum_{i=1} x_{\mu f_i}(\tilde{y}) = \frac{1}{n} \sum_{i=1} x_{\mu g_i}(\tilde{y})$ is a stationary point of $F$ in the sense of~\eqref{eq:stationary_condition}. \\
3. $y^* \in \argmin F_\mu (y)$ if and only if $ \frac{1}{n} \sum_{i=1} x_{\mu f_i}(y^*)  $ $ = \frac{1}{n} \sum_{i=1} x_{\mu g_i}(y^*) = x^* \in \argmin F(y)$.  
\end{lemma}

\noindent The minimization problem~\eqref{eq:introprob} is generally challenging due to the nonconvexity and nonsmoothness of the objective functions. In contrast, problem~\eqref{eq:intro_mu_prob} with the approximation $F_\mu$ provides an attractive surrogate: as $\nabla F_\mu$ is Lipschitz
continuous, which is a desirable property for a wide range of first-order methods. Moreover, from Lemmas~\ref{lem:moreau_properties} and~\ref{lem:problem_correspondance}, $F_\mu$ also largely preserves the geometric structure of $F$. Obtaining a stationary solution of $F_\mu$, we can recover its counterpart for $F$ via the proximal mappings $x_{\mu f_i}$ and $x_{\mu g_i}$. Thus, we focus on solving problem~\eqref{eq:intro_mu_prob}.

\subsection{Proposed $\dc$ Algorithm}\label{sec:proposed_alg}
Problem~\eqref{eq:intro_mu_prob} can be recast by creating local copies $y_i$ for all $i \in \{1,2,\dots, n\}$, of the solution $y$ to problem~\eqref{eq:intro_mu_prob} and imposing the agreement of the solutions of all the agents via
consensus constraint leading to the equivalent problem,
\begin{align}\label{eq:mu_prob_cons_form}
   \minimize_{y_1, y_2, \dots, y_n} \ &  \textstyle \frac{1}{n} \sum_{i=1}^n (\mfi(y_i) - \mgi(y_i)),\\
   \mbox{ subject to} \ & y_i = y_j, \mbox{ for all } i,j. \nonumber
\end{align} 
We develop a distributed gradient descent-based algorithm to solve problem~\eqref{eq:mu_prob_cons_form} called the $\dc$ algorithm. The algorithm proceeds in the following manner:\\
At any iteration $k$ of the algorithm, each agent $i$ maintains two estimates, an \textit{optimization variable} $y_i^{(k)}\in \mathbb{R}^p$, and a local \textit{update variable} $z_i^{(k)} \in \mathbb{R}^p$. Every iteration $k$ involves two updates: first, each agent $i$ updates $z_i^{(k)}$ via local gradient descent based on the mapping $\mFi$, with the gradient evaluated at $y_i^{(k-1)}$; next, the optimization variable $y_i^{(k)}$ is updated to an estimate which is $\eta_k$-close to the average value $\widehat{z}^{(k)}:= \frac{1}{n} \sum_{i=1}^n z_i^{(k)}$, i.e., $\|y_i^{(k)} - \widehat{z}^{(k)}\| \leq \eta^{(k)}$, using the distributed $\eta$-consensus protocol (described in detail in the next section), initialized with $z_i^{(k)}$ as the initial condition for the agent $i$ and tolerance $\eta_k$. The above algorithm updates at any agent $i$ in are summarized in the next equations: 
\begin{align}
         & \hspace{-1.86in} z_i^{(k+1)} = y_i^{(k)}  - \alpha \textstyle \nabla \frac{1}{n}\mFi(y_i^{(k)}) \label{eq:gradStep1}\\ 
               & \hspace{-1.5in} = y_i^{(k)}  -  \textstyle \frac{\alpha}{n} \nabla ( \mfi(y_i^{(k)}) - \mfi(y_i^{(k)}) ) \nonumber \\   
               &  \hspace{-1.5in} = y_i^{(k)}  - \textstyle \frac{\alpha \mu^{-1}}{n}  (x_{\mu g_i}  (y_i^{(k)}) - x_{\mu f_i}(y_i^{(k)}) ), \label{eq:gradStep2} \\   
    & \hspace{-1.85in} y_i^{(k+1)} = w_i^{(t_{k+1})}, \ \ \mbox{where} \label{eq:conStep1}\\
   \| w_i^{(t_{k+1})} - \widehat{z}^{(k+1)} \| < \eta^{(k+1)}, \ \widehat{z}^{(k+1)}& = \frac{1}{n} \sum_{i=1}^n z_i^{(k+1)} \label{eq:conStep2}
\end{align}
where, $w_i^{(t_k)}$ is the output of the $\eta$-consensus protocol and is an approximate estimate of the average $\widehat{z}^{(k)}$, $t_k$ denotes the number of communication steps utilized by the $\eta$-consensus protocol at iteration $k$ of the $\dc$ algorithm. We summarize $\dc$ in Algorithm~\ref{alg:gradcons_weakly}.

\begin{algorithm}[h]
\small
    \SetKwBlock{Input}{Input:}{}
    \SetKwBlock{Initialize}{Initialize:}{}
    \SetKwBlock{STEPONE}{STEP 1:}{}
    \SetKwBlock{STEPTWO}{STEP 2:}{}
    \SetKwBlock{Repeat}{Repeat for $ k = 0,1,2, \dots$}{}
    \newcommand{\inparallel}{\textbf{In parallel}}
    \newcommand{\until}{\textbf{until}}
    \Input{ $0 < \mu < 1/m_f$; step-size: $ 0 < \alpha \leq 1/L_{\mu F}$;\\ consensus tolerances $\{\eta^{(k)}\}_{k \geq 0}$}
    \Initialize{ 
        For each agent $i \in \mathcal{V}, y_i^{(0)} = z_i^{(0)} \in \mathbb{R}^p$;
      }
    \Repeat {
    \For{$i = 1,2,3,\dots,n $, (\inparallel)}{
        \vspace{0.05in}
        \tcc{ \small gradient descent iteration:}
        $z_i^{(k+1)} = y_i^{(k)}  - \textstyle \frac{\alpha \mu^{-1}}{n} (x_{\mu g_i}(y_i^{(k)}) - x_{\mu f_i}(y_i^{(k)}) )$ \\[1ex]
        \tcc{ \small consensus iterations: }
        $y_i^{(k+1)} \longleftarrow \eta^{(k+1)}$-consensus$(z_i^{(k+1)}, \ i \in \mathcal{V})$ \vspace{0.05in}
        }
        }
        \until \ a stopping criterion is met
        \vspace{0.01in}
        \caption{$\dc$}
        \label{alg:gradcons_weakly}
\end{algorithm}
\noindent The finite-time $\eta$-consensus protocol is described next. 

\begin{subsection}{Finite-time $\eta$-consensus Protocol}\label{sec:etaCons}
The finite-time $\eta$-consensus protocol and its variants are proposed in earlier articles \cite{yadav2007distributed, ratio_cons_lab_conf, ratio_cons_lab_tcns, melbourne2020geometry, MelbourneConsensusTac} by the authors under various practical scenarios. Here, we resort to consensus with vector-valued states \cite{melbourne2020geometry, MelbourneConsensusTac}. Consider, a set of $n$ agents connected via a directed graph $\mathcal{G}(\mathcal{V},\mathcal{E})$. The finite-time $\eta$-consensus protocol aims to design a distributed protocol so that the agents can compute an approximate estimate of the average, $\widehat{u} := \frac{1}{n} \sum_{i=1}^n u_i^{(0)}$, of their initial states $u_i^{(0)}$ in finite-time. This approximate estimate is parameterized by a small tolerance $\eta$ chosen apriori to make the estimate precise. In the $\eta$-consensus protocol, the agents maintain state variables $u_i^{(k)}  \in \mathbb{R}^p, v_i^{(k)}  \in \mathbb{R}$ and update them as: 
\begin{align}
   u_i^{(k)} &= \textstyle p_{ii} u_i^{(k-1)} +  \sum_{j\in\mathit{\mathcal{N}^-_{i}}}p_{ij}u_j^{(k-1)} \label{eq:cons_u}\\
   v_i^{(k)} &= \textstyle p_{ii}v_i^{(k-1)}+  \sum_{j\in\mathit{\mathcal{N}^-_{i}}}p_{ij}v_j^{(k-1)}\label{eq:cons_v} \\
   w_i^{(k)} &= \textstyle \frac{1}{v_i^{(k)}}u_i^{(k)}, \label{eq:cons_w}
\end{align}
where, $w_i^{(0)} = u_i^{(0)}, v_i^{(0)} = 1$ for all $i \in \mathcal{V}$ and $\mathcal{N}^{-}_i$ denotes the set of in-neighbors of agent $i$. The updates~(\ref{eq:cons_u})-(\ref{eq:cons_w}) are based on the push-sum (or ratio consensus) updates (see \cite{kempe2003gossip}). The following assumption on the graph $\mathcal{G}(\mathcal{V},\mathcal{E})$ and weight matrix $\mathbf{P} := [p_{ij}]$ is made:
\begin{assmp}\label{ass:graph_ass}
The directed graph $\mathcal{G}(\mathcal{V},\mathcal{E})$ is strongly connected and the weight matrix $\mathbf{P}$ is column-stochastic. 
\end{assmp}

\noindent Note that $\mathbf{P}$ being a column stochastic matrix allows for a distributed synthesis of the $\eta$-consensus protocol. The variable $w_i^{(k)} \in \mathbb{R}^{p}$ is an estimate of the average $\widehat{u}$ with each agent $i$ at any iteration $k$. The estimates $w_i^{(k)}$ converge to the average $\widehat{u}$ asymptotically.
\begin{theorem}(\cite{kempe2003gossip},\cite{saraswat2019distributed})\label{thm:consensusconv}
Let Assumption~\ref{ass:graph_ass} hold. Let $\{w_i^{(k)}\}_{k \geq 0}$ be generated by~(\ref{eq:cons_w}). Then $w_i^{(k)}$ asymptotically converges to $\widehat{u} = \frac{1}{n}\sum_{i=1}^n u_i^{(0)}$ for all $i \in \mathcal{V}$, i.e.,
\begin{align*}
    \textstyle \lim_{k \rightarrow \infty} w_i^{(k)} = \frac{1}{n}\sum_{i=1}^n u_i^{(0)}, \ \text{for all} \ i \in \mathcal{V}.
\end{align*}
\end{theorem}
\noindent Every agent $i \in \mathcal{V}$ maintains an additional scalar value $r_i$ to determine when the state $w_i^{(k)}$ are $\eta$-close to each other and hence from Theorem~\ref{thm:consensusconv}, $\eta$-close to $\widehat{u}$. The radius variable $w_i$ is designed to track the radius of the minimal ball that can enclose all the states $w_i^{(k)}$ (more details can be found in \cite{melbourne2020geometry, MelbourneConsensusTac}). The radius $w_i^{(k)}$, for all $i \in \mathcal{V}$, is updated as: 
\begin{align}\label{eq:radius}
    &\hspace{-0.05in} r_i^{(k)} := \nonumber \\
    &\hspace{-0.1in}  \begin{cases}
    0, & \hspace{-0.8in}\mbox{if} \  k = m\mathcal{D}, m = 0,1,2,\dots\\
     \max \limits_{j \in \mathcal{N}_i^-} \Big \{ \| w_i^{(k)} - w_i^{(k-1)}\| + r_j^{(k-1)} \Big \}, & \mbox{otherwise},  
 \end{cases}
\end{align}
where $\mathcal{D}$ is an upper bound on the diameter of the directed graph $\mathcal{G}$.
Denote, $\mathcal{B}(w_i^{(k)}, r_i^{(k)})$ as the $p$-dimensional ball of radius $r_i^{(k)}$ centered at $w_i^{(k)}$. It is established in \cite{melbourne2020geometry} that under Assumption~\ref{ass:graph_ass}, after $\mathcal{D}$ iterations of update~(\ref{eq:radius}), the ball $\mathcal{B}(w_i^{(k + \mathcal{D})}, r_i^{(k+\mathcal{D})})$ encloses the states $w_i^{(k)}$ of all the agents $i \in \mathcal{V}$. Further, it is also established that radius update sequences $\{r_i^{(k)}\}_{k \geq 0}$ converge to zero as $m \to \infty$. 
\begin{theorem}(\cite{melbourne2020geometry, MelbourneConsensusTac})\label{thm:eCons}
Let updates~\eqref{eq:cons_u}-\eqref{eq:cons_w} hold. Let $\{r_i^{(k)}\}_{k \geq 0}$ be the sequences generated by~(\ref{eq:radius}). Under Assumption~\ref{ass:graph_ass},  
\begin{align*}
  \textstyle \lim_{k \rightarrow \infty} r_i^{(k)} = 0, \ \mbox{for all} \ i \in \mathcal{V}. 
\end{align*} 
\end{theorem}

\noindent Theorem~\ref{thm:eCons} gives a criterion for termination of the consensus iterations~(\ref{eq:cons_u})-(\ref{eq:cons_w}) by utilizing the radius updates at each agent $i \in \mathcal{V}$ given by~(\ref{eq:radius}).
Using the value of the $r_i^{(k)}, k =0,1,2\dots$ each agent can ascertain the deviation from the consensus among the agents. 
As a method to detect $\eta$-consensus, at every iteration of the form $m \mathcal{D}$, for $m = 1, 2, \dots$,
$r_i^{(m\mathcal{D})}$ is compared to the tolerance $\eta$, if  $r_i^{(m\mathcal{D})} < \eta $ then all the agent states $w_i^{((m-1) \mathcal{D})}$ were $\eta$-close to $\widehat{u}$ (from Theorem~\ref{thm:consensusconv}) and the iterations~(\ref{eq:cons_u})-(\ref{eq:cons_w}) are terminated. Proposition~\ref{prop:finite-time} establishes that the $\eta$-consensus protocol converges in a finite number of iterations.


\begin{proposition}\label{prop:finite-time}
Under the Assumption~\ref{ass:graph_ass}, $\eta$-consensus is achieved in a finite number of iterations $k_\eta$ at each agent $i \in \mathcal{V}$. Moreover,  $k_\eta$ satisfies
\begin{align}
  \textstyle k_\eta \geq \ceil*{ \frac{\log (1/\eta)}{-\log (\delta_2)} +  \frac{\log (C\|\mathbf{u}^{(0)}\|_1/\delta_1)}{-\log (\delta_2)} }, \label{eq:eta_cons_num_iter}
\end{align}
where, $\mathbf{u}^{(0)} := [u_1^{(0)}, \dots, u_n^{(0)}] \in \mathbb{R}^{n \times p}$ and $C,\delta_1,\delta_2 \in (0,1)$ are constants related to the matrix $\mathbf{P}$ and graph $\mathcal{G}$.
\end{proposition}
\begin{proof}
Note, $r_i^{(k)} \rightarrow 0$ as $k \rightarrow \infty$. Thus, given $\eta > 0, i \in \mathcal{V}$ there exists finite $k_{i,\eta}$ such that for $ k \geq k_{i,\eta}, r_i^{(k)} < \eta$, for all $i \in \mathcal{V}.$ Choosing $k_\eta := \max_{1 \leq i \leq n} k_{i, \eta}$ establishes the claim. Further,~\eqref{eq:eta_cons_num_iter} is proved rigorously in \cite{gradcons_tnse}, Lemma 3.1. 
\end{proof}

\noindent To have a global detection, each agent generates a one-bit ``converged flag'' indicating its detection. The flag signal can be combined using distributed one-bit consensus updates (see \cite{melbourne2020geometry}) allowing the agents to achieve global $\eta$-consensus. 
\end{subsection}

\section{Convergence Analysis For $\dc$}\label{sec:convergence}
Let the average of the optimization variables at iteration $k$ be denoted as: $\widehat{y}^{(k)} := \frac{1}{n}\sum^{n}_{i=1} y_i^{(k)}$. Denote the gradient of the function $\mF$ evaluated at individual optimization variables of all the agents and the average $\widehat{y}^{(k)}$ as:
\begin{align}\label{eq:g_ghat_k}
  h^{(k)} &:= \textstyle \gmF(y_i^{(k)}), \ \mbox{and} \\
  \widehat{h}^{(k)} &:= \textstyle \gmF(\widehat{y}^{(k)})
\end{align}
respectively. A consequence of the $\eta$-consensus protocol is that the difference between $h^{(k)}$ and $\widehat{h}^{(k)}$ is bounded. 
\begin{lemma}\label{lem:gradbound} 
Let $\eta^{(k)}$ denote the consensus tolerance in the $\dc$ algorithm at iteration $k$, then 
\begin{align*}
 \|h^{(k)} - \widehat{h}^{(k)}\|  \leq \textstyle 2 L_{\mu F} \eta^{(k)},
\end{align*}
where $L_{\mu F}$ is the constant as defined in Lemma~\ref{lem:moreau_properties}.
\end{lemma}
\begin{proof}
Note that,
\begin{align}
    \|h^{(k)} - \widehat{h}^{(k)}\| &= \left\|\gmF(y_i^{(k)}) - \gmF(\widehat{y}^{(k)}) \right\| \nonumber \\    
    & \hspace{-0.85in} \textstyle = \frac{1}{n} \Big\|
    \left(\sum^{n}_{i=1} \nabla \mfi(y_i^{(k)}) - \sum^{n}_{i=1} \nabla \mgi(y_i^{(k)}) \right) \nonumber \\   
    & \hspace{-0.45in} \textstyle - \left( \sum^{n}_{i=1}\nabla \mfi(\widehat{y}^{(k)}) - \sum^{n}_{i=1}\nabla \mgi(\widehat{y}^{(k)}) \right)  \Big \| \nonumber\\   
    & \hspace{-0.85in} \textstyle = \frac{1}{\mu n} \Big\|
    \left(\sum^{n}_{i=1} \left( x_{\mu g_i}(y_i^{(k)}) - x_{\mu f_i}(y_i^{(k)}) \right) \right) \nonumber \\   
    & \hspace{0.3in} \textstyle - \left(\sum^{n}_{i=1} \left( x_{\mu g_i}(\widehat{y}^{(k)}) - x_{\mu f_i}(\widehat{y}^{(k)}) \right) \right)  \Big \| \nonumber\\ 
    & \hspace{-0.85in} \textstyle = \frac{1}{\mu n} \Big\|
    \left(\sum^{n}_{i=1} \left( x_{\mu f_i}(\widehat{y}^{(k)}) - x_{\mu f_i}(y_i^{(k)}) \right) \right) \nonumber \\   
    & \hspace{0.2in} \textstyle + \left(\sum^{n}_{i=1} \left( x_{\mu g_i}(y_i^{(k)}) - x_{\mu g_i}(\widehat{y}^{(k)})  \right) \right)  \Big \| \nonumber \\
    & \hspace{-0.85in} \textstyle \leq \frac{1}{n \mu(1 - \mu m_f)} \sum_{i=1}^{n} \|y_i^{(k)} - \widehat{y}^{(k)} \| + \frac{1}{n \mu} \sum_{i=1}^{n} \|y_i^{(k)} - \widehat{y}^{(k)} \| \nonumber \\
    & \hspace{-0.85in} \textstyle = \frac{2 - \mu m_f}{n(\mu - \mu^2 m_f)} \sum_{i=1}^{n} \|y_i^{(k)} - \widehat{y}^{(k)} \| \nonumber \\
    & \textstyle \hspace{-0.85in} = \frac{2 - \mu m_f}{n(\mu - \mu^2 m_f)} \sum_{i=1}^{n} \|y_i^{(k)} - \widehat{y}^{(k)} + \widehat{z}^{(k)} - \widehat{z}^{(k)} \| \nonumber \\
    & \textstyle \hspace{-0.85in} \leq \frac{L_{\mu F}}{n} \sum_{i=1}^{n} \|y_i^{(k)} - \widehat{z}^{(k)} \| + \frac{L_{\mu F}}{n} \sum_{i=1}^{n}  \| \widehat{z}^{(k)} - \widehat{y}^{(k)} \| \nonumber \\
    & \hspace{-0.85in} \textstyle \leq \frac{L_{\mu F}}{n} \sum_{i=1}^{n} (\eta^{(k)} + \eta^{(k)}) \leq 2L_{\mu F} \eta^{(k)} \nonumber,
\end{align}
where, $L_{\mu F} = \frac{2 - \mu m_f}{(\mu - \mu^2 m_f)}$ as defined in Lemma~\ref{lem:moreau_properties} and the last step follows from~(\ref{eq:conStep1}) and~(\ref{eq:conStep2}).
\end{proof}

\noindent Recall~(\ref{eq:gradStep1}), $\textstyle \frac{1}{n} \sum_{i=1}^n z_i^{(k+1)} =  \textstyle \frac{1}{n} \sum_{i=1}^n y_i^{(k)} - \alpha \nabla \mFi(y_i^{(k)}) $ $ = \widehat{y}^{(k)} - \alpha h^{(k)}$. 
Therefore, we have, $\widehat{z}^{(k+1)} = \widehat{y}^{(k)} - \textstyle \alpha h^{(k)} + \widehat{y}^{(k+1)} - \widehat{y}^{(k+1)}$
\begin{align}
   \implies \widehat{y}^{(k+1)} &= \widehat{y}^{(k)} - \textstyle \alpha h^{(k)} + a^{(k+1)}, \label{eq:closedform_upd}
\end{align}
with $ a^{(k+1)}:= \widehat{y}^{(k+1)} - \widehat{z}^{(k+1)}, \|a^{(k+1)}\| \leq \eta^{(k+1)}$. In the centralized setting, the information about the gradient of the function $\mF$, i.e. $\nabla \mF(y) = \gmF(y) = \widehat{h}_y$ is known to the solver. Here, an iteration of the (centralized) gradient descent is of the form: $\widetilde{y} = y - \alpha \widehat{h}_y$, where, $\widetilde{y}$ denotes the updated estimate of the solution. Due to Lemma~\ref{lem:gradbound}, update~(\ref{eq:closedform_upd}) can be viewed as an ``inexact'' centralized gradient descent update on the gradient $\gmF$ evaluated at the average of all the agents' optimization variables for the function $\mF$, 
\begin{align}
    \widehat{y}^{(k+1)} &= \widehat{y}(k) - \textstyle \alpha \widehat{h}^{(k)} + e^{(k)}, \ \mbox{with} \label{eq:gradstep_with_avg_grad}\\
     e^{(k)} &:= \alpha (\widehat{h}^{(k)} - h^{(k)}) + a^{(k+1)}, \label{eq:uk}\\
     \|e^{(k)}\| & \leq 2\alpha L_{\mu F} \eta^{(k)} + \eta^{(k+1)}. \label{eq:uk_norm}
\end{align} 
Therefore, the proposed $\dc$ algorithm updates are equivalent to an (approximate) centralized gradient descent update to minimize the function $F_{\mu}$.

\begin{theorem}\label{thm:weakly_convex}
Let assumptions~\ref{assm:weakly_convex} and~\ref{ass:graph_ass} hold, $\alpha \leq 1/L_{\mu F}$, and the tolerance sequence in $\dc$ is $\eta^{(k)} = 1/k^{1+\theta}$, for all $k$, with a positive constant $\theta$. Suppose $\{y_i^{(k)}, x_{\mu f_i}(y_i^{(k)}), x_{\mu g_i}(y_i^{(k)}) \}_{k \geq 0, i \in \mathcal{V}}$ be the sequences generated by  Algorithm~\ref{alg:gradcons_weakly}. Let $\widehat{y}^{(k)} := \frac{1}{n}\sum^{n}_{i=1} y_i^{(k)}$, $\widehat{\xi}^{(k)} :=  \textstyle \mu^{-1} (\frac{1}{n} $ $ \sum_{i=1}^n x_{\mu g_i}(\widehat{y}^{(k)}) - \frac{1}{n}\sum_{i=1}^n x_{\mu f_i}(\widehat{y}^{(k)}) )$, and $\xi^{(k)} :=  \textstyle \mu^{-1} (\frac{1}{n} \sum_{i=1}^n x_{\mu g_i}(y_i^{(k)}) - \frac{1}{n}\sum_{i=1}^n x_{\mu f_i}(y_i^{(k)}) )$. Then, for any positive integer K, there exists $0 \leq \overline{k} \leq K - 1$ such that 
\begin{align}
    & \textstyle \widehat{\xi}^{(\overline{k})}  \in \partial  \frac{1}{n}\sum_{i=1}^n f_i( x_{\mu f_i} (\y^{(k)}) ) -  \partial \frac{1}{n}\sum_{i=1}^n g_i( x_{\mu g_i} (\y^{(k)}) ), \nonumber \\
    &  \max \{ \| \widehat{\xi}^{(\overline{k})}\|,\| \textstyle  \frac{1}{n} \sum^{n}_{i=1} x_{\mu g_i} (\y^{(\overline{k})}) -  \frac{1}{n}\sum^{n}_{i=1} x_{\mu f_i} (\y^{(\overline{k})}) \| \}  \nonumber \\
    & \hspace{0.05in} \leq \textstyle \max \{1,\mu^{-1} \} \bigg( \frac{2 \mu^2 \left( F_\mu(\y^{(0)}) - F^* +  \textstyle \sum_{k=0}^{\infty} \frac{2 C}{k^{1+ \theta}} \right)}{\alpha K}  \bigg)^{1/2}. \nonumber
\end{align} 
\noindent Therefore, given $\varepsilon >0$, in no more than 
\begin{align}
   \hspace{-0.05in} K \leq \textstyle \bigg \lceil 
  \frac{2 \max \{1,\mu^{-1} \}^2 \mu^2 \left( F_\mu(\y^{(0)}) - F^* +  \textstyle \sum_{k=0}^{\infty} \frac{2 C}{k^{1 + \theta}} \right)}{\alpha \varepsilon^2} \bigg \rceil \hspace{-0.05in} = O(\frac{1}{\varepsilon^{2}}) \nonumber 
\end{align}
iterations, $\widehat{\xi}^{(\overline{k})}$ is an $\varepsilon$-stationary point in the sense of~\eqref{eq:epsilon_stationary_condition}, and Algorithm~\ref{alg:gradcons_weakly} estimate $\xi^{(\overline{k})}$ is $L_{\mu F} \eta^{(\overline{k})}$-close to an $\varepsilon$-stationary point in the sense of~\eqref{eq:epsilon_stationary_condition}.\\
\hspace*{0.1in} Moreover, if $F$ is level-bounded and $\dom (f_i) = \mathbb{R}^p$ for all $i$, then sequences $\{y_i^{(k)}\}$ stay bounded and for every limit point $y^\infty$ of $y_i^{(k)}$ for all $i$, we have, $\frac{1}{n}\sum_{i=1}^n x_{\mu g_i}(y^\infty) = \frac{1}{n}\sum_{i=1}^n x_{\mu f_i}(y^\infty)$ and $\frac{1}{n}\sum_{i=1}^n x_{\mu f_i}(y^\infty)$ satisfies the exact stationary condition~\eqref{eq:stationary_condition}.
\end{theorem}
\begin{proof}
\noindent For any $k \in \mathbb{Z}^{+}$, let
\begin{align}
    \widehat{\xi}^{(k)} & = \textstyle \mu^{-1} \left( \frac{1}{n}\sum^{n}_{i=1} x_{\mu g_i}(\y^{(k)}) - \frac{1}{n}\sum^{n}_{i=1} x_{\mu f_i} (\y^{(k)}) \right) \nonumber \\
    & = \textstyle \mu^{-1} \left( \y^{(k)} - \frac{1}{n}\sum^{n}_{i=1}  x_{\mu g_i}(\y^{(k)}) \right) \nonumber \\
    & \hspace{0.2in} \textstyle - \mu^{-1}\left( \y^{(k)} - \frac{1}{n}\sum^{n}_{i=1} x_{\mu f_i} (\y^{(k)}) \right). \label{eq:inequality_3}
\end{align}
From the optimality of the proximal mapping, $x_{\mu f_i} (\y^{(k)})$, 
\begin{align}
   & 0 \in \partial \{ f_i( x_{\mu f_i} (\y^{(k)}) ) + \textstyle \frac{1}{\mu} (x_{\mu f_i} (\y^{(k)}) - \y^{(k)} ) \} \nonumber \\
   & \implies \mu^{-1} ( \y^{(k)} - x_{\mu f_i} (\y^{(k)}) ) \in \partial f_i( x_{\mu f_i} (\y^{(k)}) ).  
\end{align}
Thus, 
\begin{align}
   \hspace{-0.1in}  \textstyle \frac{1}{\mu} ( \y^{(k)} - \frac{1}{n} \sum_{i=1}^n x_{\mu f_i} (\y^{(k)}) ) & \in \textstyle \partial \frac{1}{n} \sum_{i=1}^n \partial f_i( x_{\mu f_i} (\y^{(k)}) ) \nonumber \\
   & \hspace{-0.2in} \textstyle \subseteq \partial \frac{1}{n} \sum_{i=1}^n f_i( x_{\mu f_i} (\y^{(k)}) ). \label{eq:inequality_stationarity_f}
\end{align}
Similarly, from the optimality of $x_{\mu g_i} (\y^{(k)})$ 
\begin{align}
   \hspace{-0.1in} \textstyle \frac{1}{\mu} ( \y^{(k)} - \frac{1}{n} \displaystyle \sum_{i=1}^n x_{\mu g_i} (\y^{(k)}) ) \in \textstyle \partial \frac{1}{n} \sum_{i=1}^n g_i( x_{\mu g_i} (\y^{(k)}) ). \label{eq:inequality_stationarity_g}
\end{align}
Therefore, from~\eqref{eq:inequality_3},~\eqref{eq:inequality_stationarity_f} and~\eqref{eq:inequality_stationarity_g}, 
\begin{align}
   \hspace{-0.07in} \widehat{\xi}^{(k)} \in \partial \frac{1}{n} \sum_{i=1}^n f_i( x_{\mu f_i} (\y^{(k)}) ) - \partial \frac{1}{n} \sum_{i=1}^n g_i( x_{\mu g_i} (\y^{(k)}) ). \label{eq:inequality_4}
\end{align}
Let $\xi^{(k)} = \textstyle \mu^{-1} (\frac{1}{n} \sum_{i=1}^n x_{\mu g_i}(y_i^{(k)}) - \frac{1}{n}\sum_{i=1}^n x_{\mu f_i}(y_i^{(k)}) ).$
Consider, $\xi^{(k)} - \widehat{\xi}^{(k)} = \mu^{-1} \Big(\frac{1}{n}\sum_{i=1}^n x_{\mu g_i}(y_i^{(k)}) - \frac{1}{n}\sum_{i=1}^n x_{\mu g_i}(\y^{(k)}) \Big) - \mu^{-1} \big( 
    \frac{1}{n}\sum_{i=1}^n x_{\mu f_i}(y_i^{(k)}) - \frac{1}{n}\sum_{i=1}^n $ $ x_{\mu f_i}(\y^{(k)}) \big).$
Therefore, 
\begin{align}
    \|\xi^{(k)} - \widehat{\xi}^{(k)}\| &\leq \textstyle \frac{\mu^{-1}}{n}\sum_{i=1}^n \|x_{\mu g_i}(y_i^{(k)}) - x_{\mu g_i}(\y^{(k)}) \| \nonumber \\
    & \textstyle \hspace{0.2in} + \frac{\mu^{-1}}{n}\sum_{i=1}^n \| x_{\mu f_i}(y_i^{(k)}) - x_{\mu f_i}(\y^{(k)}) \| \nonumber \\
    & \leq \textstyle \frac{\mu^{-1}}{n} \Big( \sum_{i=1}^n \|y_i^{(k)} - \y^{(k)} \| \nonumber \\
    & \textstyle \hspace{0.5in} + \frac{1}{(1 - \mu m_f)} \sum_{i=1}^n \|y_i^{(k)} - \y^{(k)} \| \Big) \nonumber \\
    & \textstyle \leq \frac{L_{\mu F}}{n} \sum_{i=1}^n \eta^{(k)} = L_{\mu F} \eta^{(k)}, \label{eq:xi-xihat-close}
\end{align}
where we used Lipschitz continuity properties of the proximal mappings listed in Lemma~\ref{lem:moreau_properties}. 

Since, $\mFi$ is $L_{\mu F}$-Lipschitz differentiable. It can be seen that $F_\mu = \mF$ is Lipschitz differentiable with constant $\frac{1}{n} \sum_{i=1}^n L_{\mu F} = L_{\mu F}$. With $\alpha \leq 1/L_{\mu F}$ and~\eqref{eq:gradstep_with_avg_grad},
    \begin{align*}
       & F_\mu(\y^{(k+1)}) \\
       &  \leq F_\mu(\y^{(k)}) + \langle \widehat{h}^{(k)}, \y^{(k+1)} - \y^{(k)} \rangle  + \textstyle \frac{L_{\mu F}}{2} \|\y^{(k+1)} - \y^{(k)}\|^2 \\
       & = F_\mu(\y^{(k)}) - \textstyle \frac{1}{\alpha} \langle \y^{(k+1)} - \y^{(k)} - e^{(k)}, \y^{(k+1)} - \y^{(k)} \rangle  \\
       & \hspace{0.2in} + \textstyle \frac{L_{\mu F}}{2} \|\y^{(k+1)} - \y^{(k)}\|^2 \\
       & = F_\mu(\y^{(k)}) + \textstyle \big(\frac{L_{\mu F}}{2} - \frac{1}{\alpha} \big)\|\y^{(k+1)} - \y^{(k)}\|^2 \\
       & \hspace{0.2in} + \textstyle \frac{1}{\alpha} \langle  e^{(k)}, \y^{(k+1)} - \y^{(k)} \rangle.
    \end{align*}
This implies,
\begin{align*}
    & \textstyle \big( \frac{1}{\alpha} - \frac{L_{\mu F}}{2} \big)\|\y^{(k+1)} - \y^{(k)}\|^2 \\
    & \leq F_\mu(\y^{(k)}) - F_\mu(\y^{(k+1)}) + \textstyle  \frac{1}{\alpha} \| e^{(k)} \| \| \y^{(k+1)} - \y^{(k)} \| \\
    & \leq F_\mu(\y^{(k)}) - F_\mu(\y^{(k+1)}) + \textstyle  \frac{1}{\alpha} \| e^{(k)} \| \| e^{(k)} - \alpha \widehat{h}^{(k)} \| \\
    & \leq F_\mu(\y^{(k)}) - F_\mu(\y^{(k+1)}) + \textstyle  \frac{1}{\alpha} \| e^{(k)} \|^2 + \| e^{(k)}\| \|\widehat{h}^{(k)} \| \\
    & \leq F_\mu(\y^{(k)}) - F_\mu(\y^{(k+1)}) + \textstyle  C_1 (\eta^{(k)})^2 + C_2 \eta^{(k)},\\
    & \leq F_\mu(\y^{(k)}) - F_\mu(\y^{(k+1)}) + \textstyle  C \eta^{(k)}
\end{align*}
where $C_1 := \frac{1}{\alpha}(2 \alpha L_{\mu F} + 1)^2$, $C_2 := \sup_k \| \widehat{h}^{(k)} \|(2 \alpha L_{\mu F} + 1)$, $C := C_1 + C_2$, and we used~\eqref{eq:gradstep_with_avg_grad}-\eqref{eq:uk_norm}. Therefore, 
\begin{align*}
    & F_\mu(\y^{(k)}) - F_\mu(\y^{(k+1)}) + \textstyle  C \eta^{(k)}\\
    & \geq \textstyle  \big( \frac{1}{\alpha} - \frac{L_{\mu F}}{2} \big)\|\y^{(k+1)} - \y^{(k)}\|^2 \\
    & \geq \textstyle \frac{1}{2\alpha} \|\y^{(k+1)} - \y^{(k)}\|^2 \geq \textstyle \frac{\alpha}{2} \|\widehat{h}^{(k)}\|^2 - \frac{1}{2} C_1 (\eta^{(k)})^2.
\end{align*}
Thus, 
\begin{align*}
    & \hspace{-0.25in} F_\mu(\y^{(k)}) - F_\mu(\y^{(k+1)}) + \textstyle  \frac{3}{2} C \eta^{(k)} \geq \textstyle \frac{\alpha}{2} \|\widehat{h}(k)\|^2 \\
    & \hspace{-0.25in} = \textstyle \frac{\alpha}{2 \mu^2} \| \frac{1}{n} \sum^{n}_{i=1} x_{\mu g_i} (\y^{(k)}) - \frac{1}{n} \sum^{n}_{i=1} x_{\mu f_i} (\y^{(k)}) \|^2.
\end{align*}
Summing the above inequality over $k = 0,1,\dots, K-1$, for some positive integer $K-1$, gives,
\begin{align}
    & \hspace{-0.1in} \textstyle \sum_{k=0}^{K-1} \| \frac{1}{n} \sum^{n}_{i=1} x_{\mu g_i} (\y^{(k)}) - \frac{1}{n} \sum^{n}_{i=1} x_{\mu f_i} (\y^{(k)}) \|^2 \nonumber \\    
    & \hspace{-0.1in} \leq \textstyle \frac{2 \mu^2}{\alpha} \left( F_\mu(\y^{(0)}) - F_\mu(\y^{(K)}) +  \textstyle \frac{3C}{2}  \sum_{k=0}^{K-1} \eta^{(k)} \right) \nonumber\\
    & \hspace{-0.1in} \leq  \textstyle \frac{2 \mu^2}{\alpha} \left( F_\mu(\y^{(0)}) - F^* +  \textstyle \frac{3C}{2}  \sum_{k=0}^{K-1} \eta^{(k)} \right) \nonumber \\
    & \hspace{-0.1in} \leq  \textstyle \frac{2 \mu^2}{\alpha} \left( F_\mu(\y^{(0)}) - F^* +  \textstyle \sum_{k=0}^{\infty} \frac{2 C}{k^{1+ \theta}} \right) .
    \label{eq:inequality_1}
\end{align}
Let $\overline{k} := \argmin \limits_{k = 0,1,\dots, K}  \Big \| \frac{1}{n} \sum^{n}_{i=1} x_{\mu g_i} (\y^{(k)}) -  \frac{1}{n} \sum^{n}_{i=1} $ $ x_{\mu f_i} (\y^{(k)}) \Big \|^2$. Therefore, from~\eqref{eq:inequality_1},
\begin{align}
    & \textstyle \left \| \frac{1}{n} \sum^{n}_{i=1} x_{\mu g_i} (\y^{(\overline{k})}) - \frac{1}{n} \sum^{n}_{i=1} x_{\mu f_i} (\y^{(\overline{k})}) \right \| \nonumber \\ 
    & \textstyle \leq \bigg( \frac{2 \mu^2 \left( F_\mu(\y^{(0)}) - F^* +  \textstyle \sum_{k=0}^{\infty} \frac{2 C}{k^{1+ \theta}} \right)}{\alpha K}  \bigg)^{1/2}. \label{eq:inequality_2}
\end{align}
\noindent Thus, from~\eqref{eq:inequality_3},~\eqref{eq:inequality_4}, and~\eqref{eq:inequality_2} with the claimed upper on bound $K$ in Theorem~\ref{thm:weakly_convex} we get,
\begin{align}
    & \hspace{-0.12in} \max \{ \| \widehat{\xi}^{(\overline{k})} \| ,\| \textstyle  \frac{1}{n} \sum^{n}_{i=1} x_{\mu g_i} (\y^{(\overline{k})}) -  \frac{1}{n}\sum^{n}_{i=1} x_{\mu f_i} (\y^{(\overline{k})}) \| \}  \nonumber \\
    & \hspace{-0.15in} \leq \textstyle \max \{1,\mu^{-1} \} \bigg( \frac{2 \mu^2 \left( F_\mu(\y^{(0)}) - F^* +  \textstyle \sum_{k=0}^{\infty} \frac{2 C}{k^{1+ \theta}} \right)}{\alpha K}  \bigg)^{1/2} \leq \varepsilon. \nonumber
\end{align}

Since $F$ is level-bounded and $\dom (f_i) = \mathbb{R}^p$ for all $i$, then from \cite{sun2023algorithms}, Proposition 5, $F_\mu$ is level-bounded. As $\{ F_\mu (\y^{(k)}) \}_{k \in \mathbb{N}}$ is monotonically non-increasing, the sequence, $\{\y^{(k)} \}_{k \in \mathbb{N}}$ is bounded and therefore, has at least one limit point $\y^\infty$. Let $\{\y^{(k_j)} \}_{j \in \mathbb{N}}$ denote the subsequence converging to $\y^\infty$. Since, $x_{\mu f_i}$ and $x_{\mu g_i}$ are continuous,~\eqref{eq:inequality_1} implies $ \frac{1}{n}\sum_{i=1}^n x_{\mu f_i}(\y^\infty) =  \frac{1}{n}\sum_{i=1}^n x_{\mu g_i}(\y^\infty)$. As $g_i$ are continuous, $\lim_{j \to \infty}  \frac{1}{n} \sum_{i=1}^n g_i(x_{\mu g_i}(\y^{(k_j)})) =  \frac{1}{n}\sum_{i=1}^n g_i(x_{\mu g_i}(\y^\infty))$; in addition, 
\begin{align*}
    \textstyle  \frac{1}{n}\sum_{i=1}^n f_i(x_{\mu f_i}(\y^\infty)) & \textstyle \leq \liminf \limits_{j \to \infty}  \frac{1}{n}\sum_{i=1}^n f_i(x_{\mu f_i}(\y^{(k_j)})) \\
    & \hspace{-0.25in} \textstyle \leq \limsup_{j \to \infty}  \frac{1}{n}\sum_{i=1}^n f_i(x_{\mu f_i}(\y^{(k_j)})) \\
    & \hspace{-1.1in} \textstyle  \leq \lim_{j \to \infty}  \frac{1}{n}\sum_{i=1}^n f_i(x_{\mu f_i}(\y^\infty)) \\
    & \hspace{-0.05in} \textstyle  + \frac{1}{2 \mu n} \sum_{i=1}^n \|x_{\mu f_i}(\y^\infty) - \y^{(k_j)} \|^2 \nonumber \\
    & \hspace{-0.05in} \textstyle  - \frac{1}{2 \mu n}  \sum_{i=1}^n \|x_{\mu f_i}(\y^{(k_j)}) - \y^{(k_j)} \|^2 \\
    & \hspace{-1.1in}  = \textstyle  \frac{1}{n}\sum_{i=1}^n f_i(x_{\mu f_i}(\y^\infty)),
\end{align*}
where the first inequality is due to the lower-semicontinuity of $\frac{1}{n}\sum_{i=1}^n f_i$ and the last inequality is due to the optimality of $x_{\mu f_i}(\y^\infty)$ in each Moreau envelope evaluation, and therefore we also have $\lim_{j \to \infty} \frac{1}{n} \sum_{i=1}^n f_i (x_{\mu f_i}(\y^{(k_j)})) = \lim_{j \to \infty} \frac{1}{n}\sum_{i=1}^n f_i (x_{\mu f_i}(\y^\infty))$. 
Due to $\eta$-consensus there exists $y^\infty$ such that $\lim_{j \to \infty} y_i^{(k_j)} = y^\infty$ for all $i = 1,\dots, n$. Further, for all $i$, $\| x_{\mu f_i} (y_i^{(k_j)}) - x_{\mu f_i} (\y^{(k_j)}) \| \leq \frac{1}{1 - \mu m_f} \| y_i^{(k_j)} - \y^{(k_j)} \| \leq \frac{1}{1 - \mu m_f} \eta^{(k_j)}$. Thus, $x_{\mu f_i} (y_i^{(k_j)}) \to x_{\mu f_i} (\y^\infty)$ for all $i = 1,2,\dots,n$. Since, $g_i$ are continuous, $\lim_{j \to \infty} g_i(x_{\mu g_i}(y_i^{(k_j)})) = g_i(x_{\mu g_i}(\y^\infty))$ for all $i = 1,2,\dots, n$. Thus, $\lim_{j \to \infty} \frac{1}{n} \sum_{i=1}^n g_i(x_{\mu g_i}(y_i^{(k_j)})) = \frac{1}{n}\sum_{i=1}^n g_i(x_{\mu g_i}(\y^\infty))$. Finally, taking the limit on~\eqref{eq:inequality_4} and~\eqref{eq:xi-xihat-close} along the appropriate subsequences gives $\xi$ satisfy~\eqref{eq:stationary_condition} in the limit.
\end{proof}


\begin{section}{Numerical Simulations}\label{sec:simulation}
This section presents a simulation study for the proposed $\dc$ algorithm. We consider a network of $10$ agents with the underlying communication network generated using the Erdos-Renyi model \cite{erdHos1960evolution} with a connectivity probability of $0.2$. We consider a $\ell_{1-2}$ regularized distributed least squares problem \cite{ahn2017difference}:
\begin{align*}
    \minimize_{x \in \mathbb{R}^p} \textstyle \ \frac{1}{10} \sum_{i=1}^{10} \frac{1}{2} \|A_i x - b_i\|^2 + \rho \|x\|_1 - \rho \|x\|_2.
\end{align*}
The problem data is generated as follows: first, we create the data matrices $A_i \in \mathbb{R}^{m \times p}$, for all $i$, with entries generated from the standard Gaussian distribution. Then the columns of the matrices are normalized to the unit norm. We generate $x^* \in \mathbb{R}^p$ such that $x^*$ is a sparse vector with $s$ non-zero entries. The non-zero entries of $x^*$ are drawn from the standard Gaussian distribution. Finally, the vectors $b_i \in \mathbb{R}^m$ are determined as $b_i = A_i x^* + 0.01 \zeta_i$, where, $\zeta_i \in \mathbb{R}^m$ has standard Gaussian entries. $F^* := \frac{1}{10} \sum_{i=1}^{10} \frac{1}{2} \|A_i x^* - b_i\|^2 + \rho \|x^*\|_1 - \rho \|x^*\|_2$. The parameters used during the simulations are: $m = 720, p = 2560, \rho = 0.1, \alpha = 0.01, \mu = 1/\lambda_{\max}(A_i^\top A_i)$, where $\lambda_{\max}(M)$ is the maximum eigenvalue of matrix $M$. We choose $\eta^{(k)} = 1/k^{1.1}$ as the tolerance sequence in Algorithm~\ref{alg:gradcons_weakly}.

\begin{figure}[h]
\centering
      \includegraphics[scale=0.25,trim={0cm 4.6cm 0.1cm 5cm},clip]{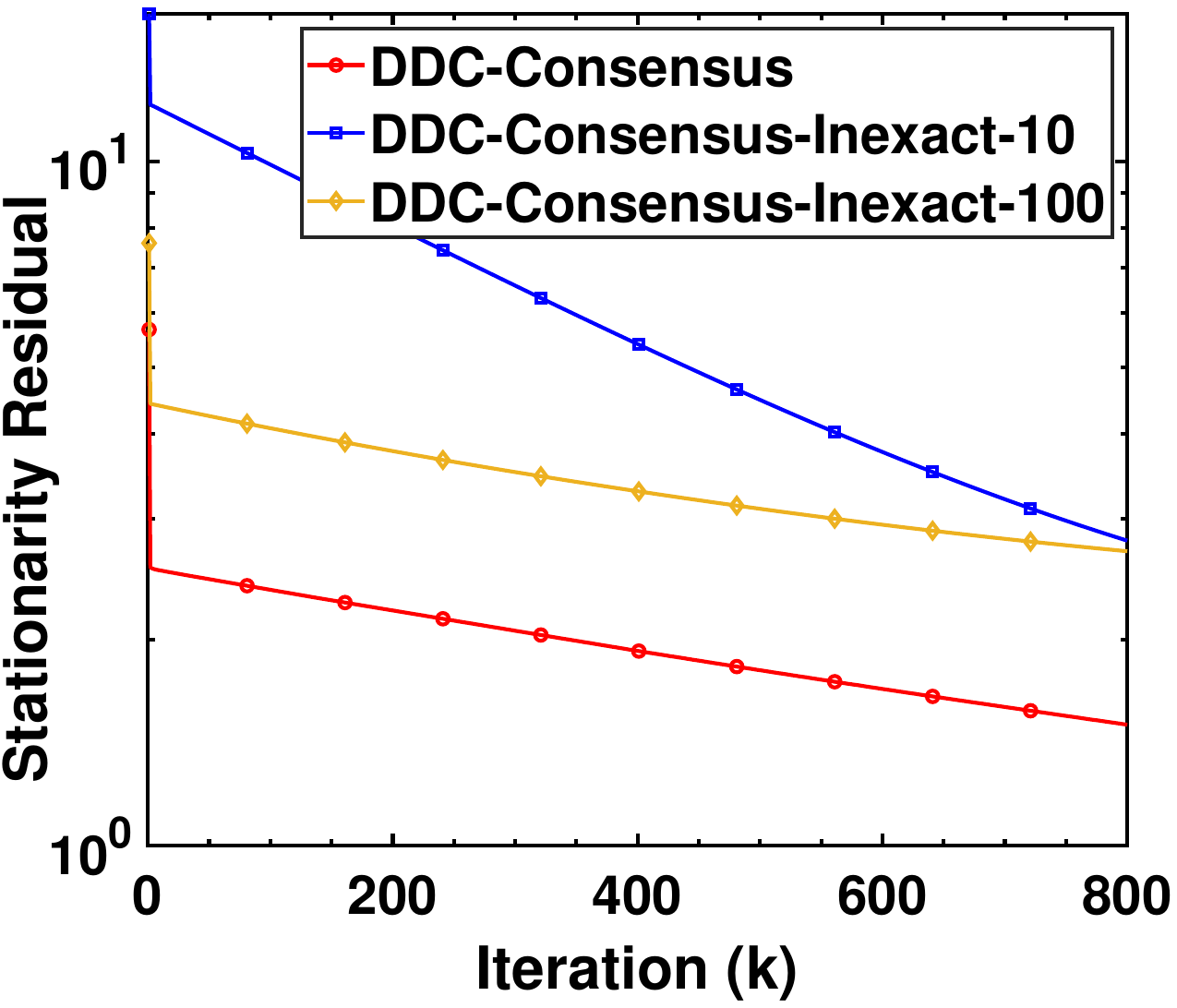}
      \caption{Stationarity residual comparison: proposed $\dc$ and $\dci$-$q$}
      \label{fig:station_inexact_compare}
\end{figure}

We compare the performance of the proposed algorithm $\dc$ with its following two relaxations: (i) $\dci$-$q$, where we solve the minimization problems~\eqref{eq:individual_moreau_f} and~\eqref{eq:individual_moreau_g} inexactly via a gradient descent method for $q$ iterations, (ii) $\dcm$, where the $\eta$-consensus is replaced with one step mixing (applying one iteration of the updates~\eqref{eq:cons_u}-\eqref{eq:cons_w} to $y_i$ estimates of the optimization algorithms). To demonstrate the performance of the proposed algorithm we present utilize the following residual metrics: \textit{Solution Residual}$(k)$: $\frac{1}{10}  \sum_{i=1}^{10} \| y_i^{(k)} - x^* \|_2$, \textit{Stationarity Residual}$(k)$: $\frac{1}{10} \sum_{i=1}^{10} \| x_{\mu g_i}(y_i^{(k)}) - x_{\mu f_i}(y_i^{(k)}) \|_2$, \textit{Objective Residual}$(k)$: $F(y_i^{(k)}) - F^*$, \textit{Consensus Residual}$(k)$: $\frac{1}{10} \sum_{i=1}^n \sum_{j=1}^n \| y_i^{(k)} - y_j^{(k)} \|_2$.

\subsection{Comparison with  $\dci$-$q$}
Figs.~\ref{fig:station_inexact_compare}-\ref{fig:sol_inexact_compare} demonstrate the performance comparisons between the proposed $\dc$ algorithm and two instances of $\dci$-$q$ that are denoted as $\dci$-$10$ and $\dci$-$100$ respectively. This comparison aims to determine a trade-off between the accuracy and computational burden of the proximal minimization in the proposed $\dc$ algorithm. As expected the $\dc$ algorithm performs the best with respect to all the performance metrics. Fig.~\ref{fig:station_inexact_compare} shows the stationarity metric that gives an estimate of the rate of convergence in the sense of Theorem~\ref{thm:weakly_convex}. 

\begin{figure}[t]
\centering
      \includegraphics[scale=0.25,trim={0cm 4.6cm 0cm 4.4cm},clip]{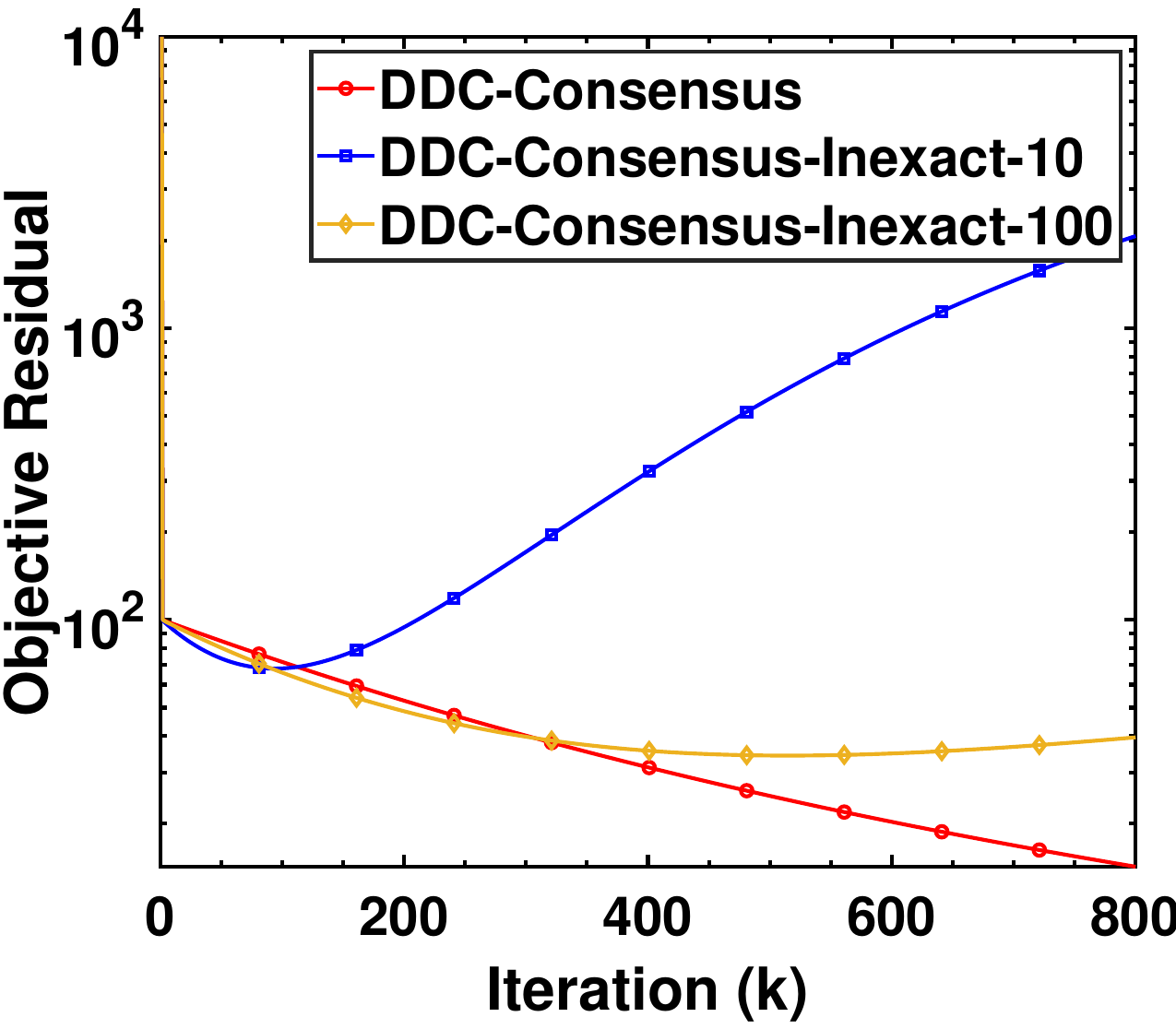}
      \caption{Objective residual comparison: proposed $\dc$ and $\dci$-$q$}
      \label{fig:fres_inexact_compare}
\end{figure}

\begin{figure}[h]
\centering
      \includegraphics[scale=0.25,trim={0.1cm 4.6cm 0cm 4.4cm},clip]{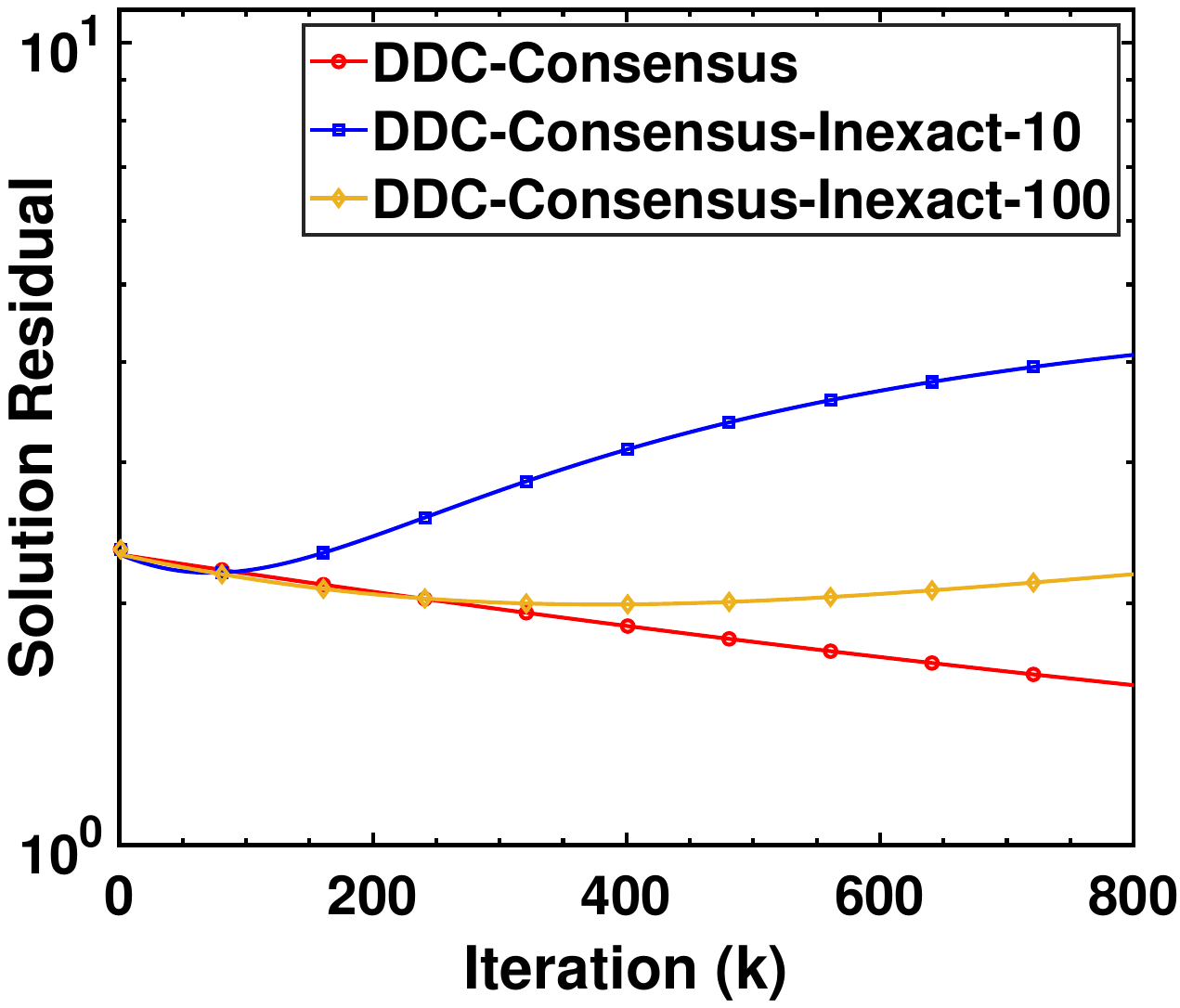}
      \caption{Solution residual comparison: proposed $\dc$ and $\dci$-$q$}
      \label{fig:sol_inexact_compare}
\end{figure}
\begin{figure}[b]
\centering
      \includegraphics[scale=0.25,trim={0.1cm 5.2cm 0.4cm 5cm},clip]{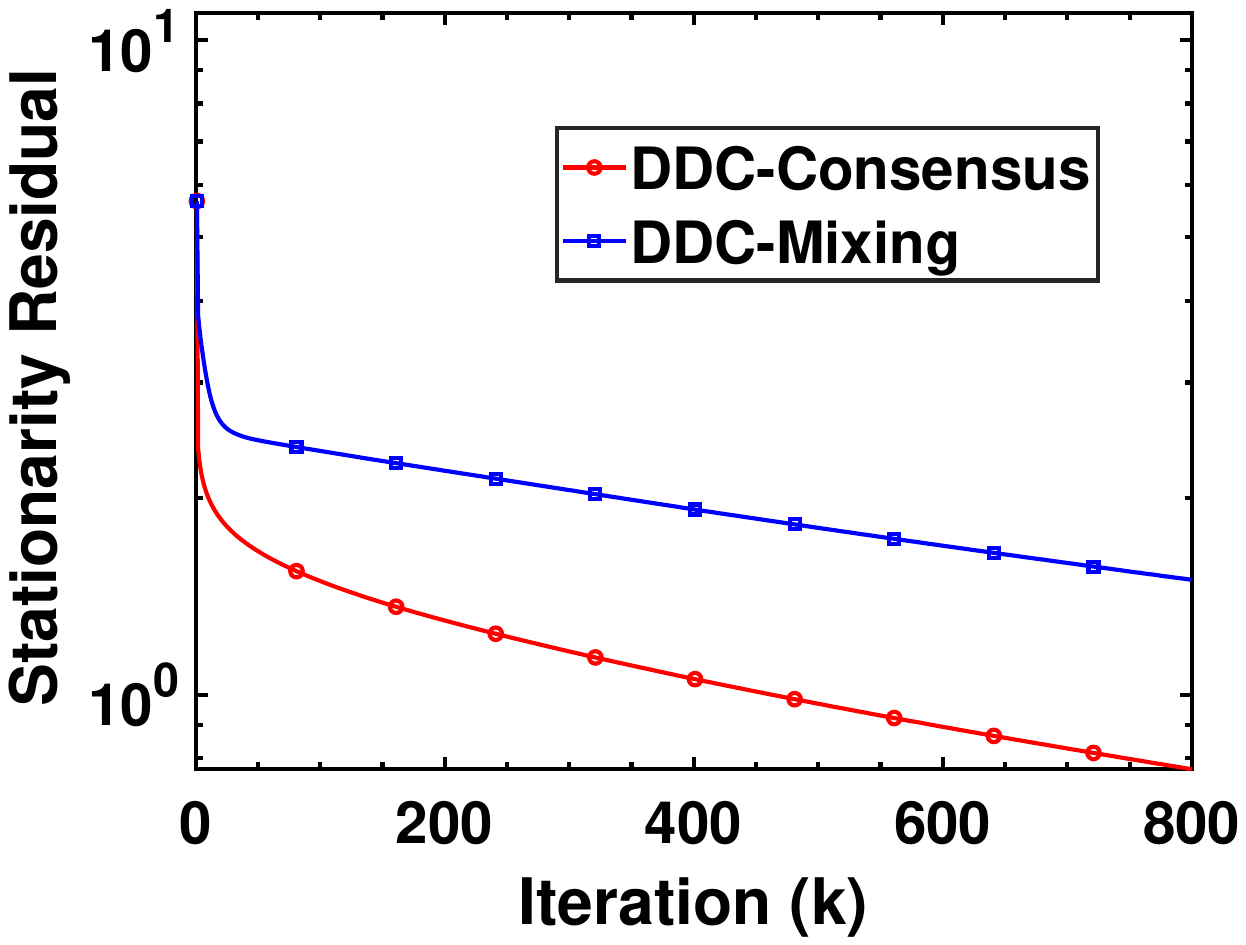}
      \caption{Stationarity residual: $\dc$ and $\dcm$}
      \label{fig:station_cons_compare}
\end{figure}
Figs.~\ref{fig:fres_inexact_compare} and~\ref{fig:sol_inexact_compare} show that gradient descent with a low computational footprint can perform relatively well compared to the exact minimization in the $\dc$ algorithm. However, it is also established that if the proximal minimization is not performed with good accuracy the algorithm iterates to diverge as in the case of $\dci$-$10$. 

\subsection{Comparison with  $\dcm$}
Figs.~\ref{fig:station_cons_compare}-\ref{fig:cons_res_cons_compare} present comparison plots showing the different metrics for both the $\dc$ and $\dcm$ algorithms. 

\begin{figure}[h]
\centering
      \includegraphics[scale=0.25,trim={0.1cm 5.2cm 0.3cm 5cm},clip]{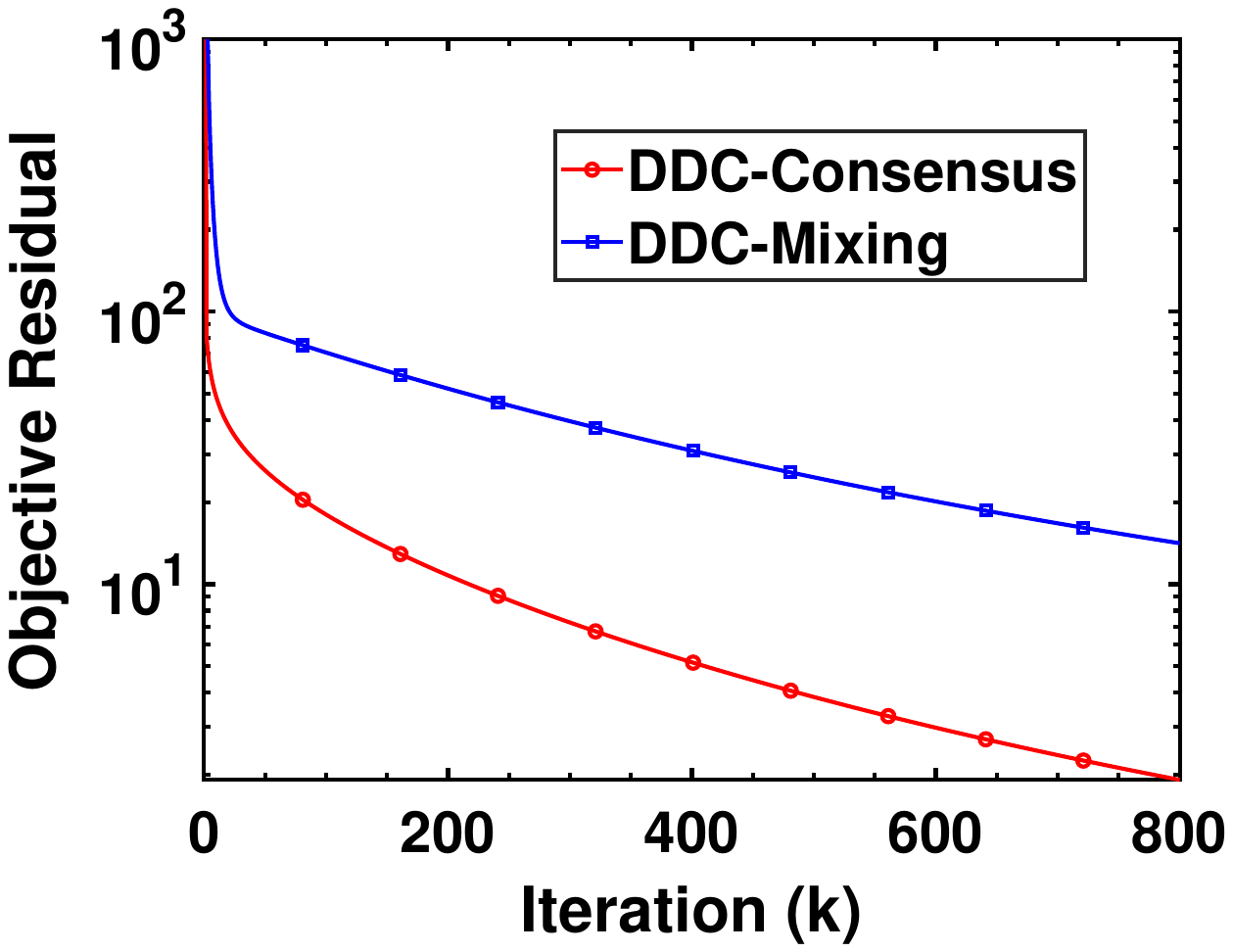}
      \caption{Objective residual: $\dc$ and $\dcm$}
      \label{fig:fres_cons_compare}
\end{figure}

The simulation study in this section aims to study the effect of the agreement quality between the agent estimates in the $\dc$ algorithm on its performance. The $\dcm$ algorithm uses \textit{only} one information exchange step among the agents and is theoretically the lower bound on information aggregation. Note, that both $\dc$ and $\dcm$ algorithms perform exact proximal minimization steps. 

\begin{figure}[h]
\centering
      \includegraphics[scale=0.25,trim={0.1cm 5.2cm 0.3cm 5cm},clip]{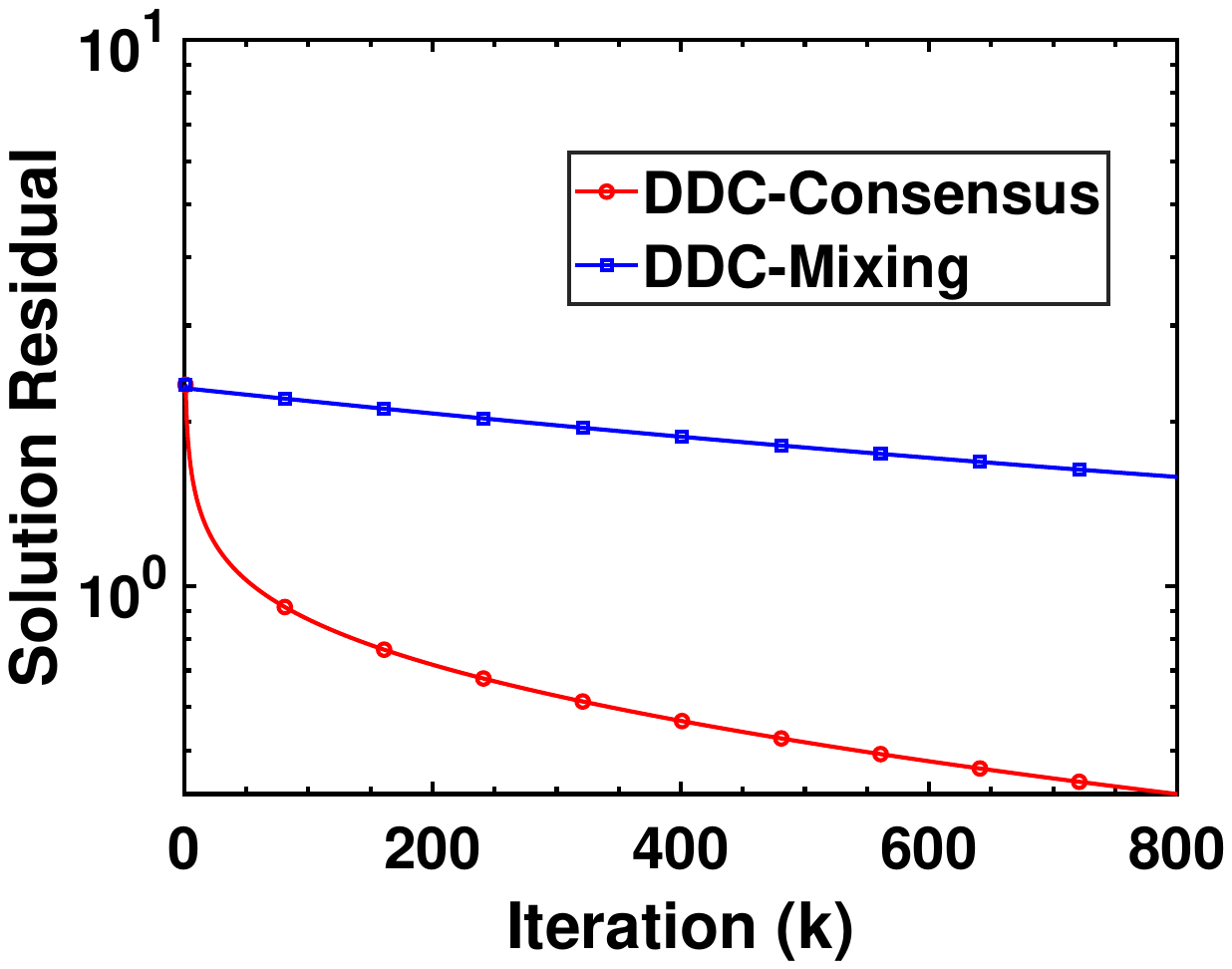}
      \caption{Solution residual comparison: proposed $\dc$ and $\dcm$}
      \label{fig:sol_cons_compare}
\end{figure}

Figs.~\ref{fig:station_cons_compare}-\ref{fig:sol_cons_compare} show that the $\eta$-consensus protocol allows for a higher degree of agreement between the agent estimates during the algorithm iterations lead to better convergence properties of the proposed $\dc$ algorithm compared to the $\dcm$ algorithm. Furthermore, the $\eta$-consensus protocol adaptively controls the agent mismatch and leads to the consensus constraint satisfaction at every iterate of the proposed $\dc$ algorithm as seen in Fig.~\ref{fig:cons_res_cons_compare}.
\begin{figure}[b]
\centering
      \includegraphics[scale=0.25,trim={0.4cm 5.2cm 0.27cm 5cm},clip]{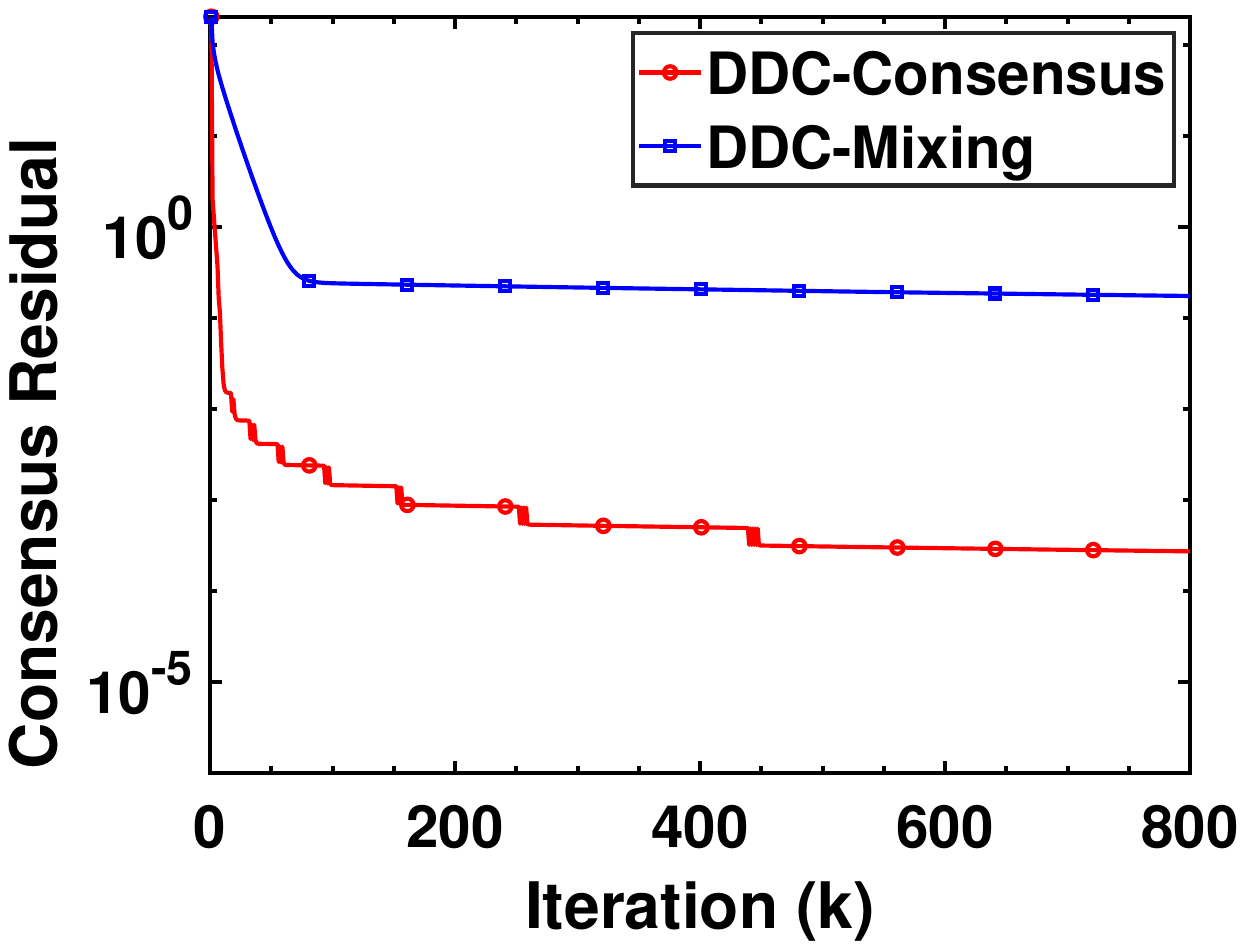}
      \caption{Consensus residual: $\dc$ and $\dcm$}
      \label{fig:cons_res_cons_compare}
\end{figure}

\end{section}

\section{Conclusion}\label{sec:conclusion}
We focused on a distributed DC optimization problem with local objective functions given by the difference between a nonsmooth weakly convex function and a nonsmooth convex function. Based on the gradient of the smooth approximations of the weakly convex and convex functions we developed the distributed $\dc$ algorithm. In the $\dc$ algorithm each agent performs a local gradient descent step and updates its estimates using a finite-time $\eta$-consensus protocol. The agent estimates are shared over a directed graph and updated via a column stochastic matrix; giving $\dc$ the desirable properties of allowing for non-symmetric communication and distributed synthesis. We established that the $\dc$ algorithm converges to a stationary point of the nonconvex distributed optimization problem. The numerical performance of the $\dc$ algorithm is presented in a simulation study solving the DC-regularized least squares problem. Results showing the behavior of $\dc$ and its relaxations $\dci$-$q$ and $\dcm$ concerning the residual metrics demonstrate the efficacy of the proposed $\dc$ algorithm.

\bibliography{references}
\end{document}